\documentclass[preprint]{elsarticle}
\usepackage{mathrsfs}
\usepackage{bm}\makeatletter \oddsidemargin -.1in \evensidemargin -.1in
\newcommand{\singlespacing}{\let\CS=\@currsize\renewcommand{\baselinestretch}{1}\tiny\CS}
\newcommand{\doublespacing}{\let\CS=\@currsize\renewcommand{\baselinestretch}{1.27}\tiny\CS}
\makeatletter
\newcommand{\setword}[2]{%
	\phantomsection
	#1\def\@currentlabel{\unexpanded{#1}}\label{#2}%
}
\makeatother

\newtheorem{definition}{Definition}[section]

\newtheorem{theorem}{Theorem}[section]
\newtheorem{proposition}{Proposition}[section]

\newproof{proof}{Proof}
\newtheorem{lemma}{Lemma}[section]

\newtheorem{remark}{Remark}[section]

\usepackage{xcolor}
\usepackage{microtype}
\usepackage{times}
\usepackage{graphicx,epsfig}
\usepackage{amsmath}
\usepackage{amssymb}
\usepackage{enumitem}
\usepackage{mathtools}
\usepackage{float}
\everymath{\displaystyle}
\usepackage[linkcolor=blue, urlcolor=blue, citecolor=blue,
colorlinks, bookmarks]{hyperref}
\date{}

\numberwithin{equation}{section}
\journal{Journal...}

\begin{document}
	
	\begin{frontmatter}
		
		\title {\textbf{\Large Recurrent Histopolating Fractal Functions}}

		 \author[1]{Pavjeet Singh}
		 \ead{pavjeetsingh23@gmail.com}
		 	\author[1]{Vinay Kumar}
		 \ead{Vinaykp20601@gmail.com}
	\author[1]{S.K. Katiyar\corref{mycorrespondingauthor}}
		\cortext[mycorrespondingauthor]{Corresponding author}
		\ead{katiyarsk@nitj.ac.in}
	
		 \address[1]{Department of Mathematics and Computing, DR BR Ambedkar National Institute of Technology, Jalandhar, Punjab, India, 144008}

		\begin{abstract}
		Recurrent fractal interpolation extends the usual fractal interpolation framework by introducing a recurrence structure. Fractal histopolation, on the other hand, is concerned with the area-matching properties of fractal functions. In this paper, we combine these two approaches and introduce a recurrent fractal histopolants.

			
			
		\end{abstract}

	\end{frontmatter}
    
		\section{Introduction}
	
	Fractal interpolation provides a recursive alternative to classical
	interpolation methods for the approximation and representation of
	irregular data. The concept of a fractal interpolation function (FIF)
	was introduced by Barnsley \cite{barnsley1986fractal} using iterated
	function systems (IFSs). Fractal interpolation has developed in several
	directions. Smoothness and calculus were studied
	in \cite{barnsley1989calculus}, while fractal polynomial interpolation  \cite{navascues2005fractal}. Function-space constructions and approximation properties of fractal functions have also been investigated in \cite{viswanathan2016associate,massopust2014fractal}. Further generalizations include perturbation error analysis of FIFs, FIFs with variable parameters, and multivariate FIFs \cite{wang2008,wang2013fractal,pandey2026generalized}.

	Recurrent iterated function systems provide an important extension of
	the usual IFS framework. Barnsley, Elton, and Hardin
	\cite{barnsley1989recurrent} introduced recurrent IFSs in which the maps are governed by a finite-state recurrence structure. In the interpolation setting, this
	framework leads to recurrent fractal interpolation functions (RFIFs),
	whose graphs arise as invariant sets of suitable recurrent iterated
	function systems. Several aspects of RFIFs and their generalizations
	have subsequently been studied. Constructions of recurrent fractal
	interpolation functions and surfaces \cite{bouboulis2006construction,liang2021construction}, while existence
	and fractal dimension of general RFIFs \cite{ruan2021existence}. Smoothness, stability,
	and further extensions have also been discussed in
	\cite{ri2021smoothness,attia2023construction,ro2024analytical}.
	
	Interpolation is based on the reproduction of prescribed function
	values at a finite set of knots. A related approximation problem is
	\emph{histopolation}, in which integral data, rather than pointwise
	values, are prescribed. Classical forms of histopolation, particularly
	spline-based histopolation, have been studied in
	\cite{fischer2007comonotone,siewer2008histopolating}. Discontinuous
	fractal functions investigated in \cite{navascues2014fractal}, showing
	that fractal constructions need not always be restricted to continuous
	functions. More recently, Barnsley and Viswanathan
	\cite{barnsley2023histopolating} introduced fractal histopolation by combining area-matching conditions with the self-referential structure of fractal functions. Their framework allows bounded integrable fractal functions that need not be continuous or interpolatory
	
	

	\section{Recurrent Continuous Fractal Interpolation Function: Revisited} \label{sec2}
	
	In this section, we briefly revisit the construction of RFIFs. Let $\Delta=\{(x_i,y_i):i=0,1,\ldots,N\}$
	be a data set with $x_0<x_1<\cdots<x_N,$
	and let $	I=[x_0,x_N].$ For each $i=1,2,\ldots,N$, define $I_i=[x_{i-1},x_i]$
	and $J_i=[x_{j(i)},x_{l(i)}],$ where $	0\leq j(i)<l(i)\leq N,$
	such that $	x_i-x_{i-1}<x_{l(i)}-x_{j(i)}.$ Let $	L_i:J_i\longrightarrow I_i$ satisfying $	L_i(x_{j(i)})=x_{i-1}, L_i(x_{l(i)})=x_i,$
	Thus, the length of the domain interval $J_i$ is strictly larger than interval $I_i$. Let  $K=J_i\times \mathbb{R}$ For each $i=1,2,\ldots,N$, consider $W_i:K\longrightarrow I_i \times \mathbb{R}$ defined as $W_i(x,y)=\bigl(L_i(x),F_i(x,y)\bigr),$ where $	L_i(x)=a_i x+b_i$ and $F_i:K \rightarrow \mathbb{R}$ defined as $$F_i(x,y)=c_i x+\alpha_i y+e_i$$ such that $	|\alpha_i|<1$
	and $	W_i(x_{j(i)},y_{j(i)})
	=
	(x_{i-1},y_{i-1}),$
	and $	W_i(x_{l(i)},y_{l(i)})
	=
	(x_i,y_i).$
    So, the constructed recurrent iterated function system (RIFS) is;
	\begin{equation}\label{RIFS}
		\left\{ K;W_i;i=1,2,\cdots,N \right\}.
	\end{equation}
    Since $L_i$ is affine, we have $a_i
	=
	\frac{x_i-x_{i-1}}
	{x_{l(i)}-x_{j(i)}}.$
	Therefore, $0<a_i<1.$ The recurrent dependence among the interpolation intervals is described by
	the connection matrix $C=(C_{im})_{N\times N},$ where
	\[
	C_{im}
	=
	\begin{cases}
		1, & I_m\subseteq J_i,\\[4pt]
		0, & \text{otherwise}.
	\end{cases}
	\]
	For each $i=1,2,\ldots,N$, define $\mathcal{I}(i)
	=
	\{m\in\{1,2,\ldots,N\}:C_{im}=1\}.$ Thus, $J_i
	=
	\bigcup_{m\in\mathcal{I}(i)}I_m.$
	Let $\mathcal{H}(\mathcal{K})$ denote the family of all nonempty compact subsets of
	$\mathcal{K}$, endowed with the Hausdorff metric. Consider the product space $\mathcal{H}(I\times \mathbb{R})^N.$ For $	A=(A_1,A_2,\ldots,A_N)\in\mathcal{H}(I\times \mathbb{R})^N,$ define the operator $\mathcal{W}:\mathcal{H}(I\times \mathbb{R})^N \longrightarrow \mathcal{H}(I\times \mathbb{R})^N$
	by
	\[
	\mathcal{W}(A_1,A_2,\ldots,A_N)
	=
	\left(
	\bigcup_{m\in\mathcal{I}(1)}W_1(A_m),
	\ldots,
	\bigcup_{m\in\mathcal{I}(N)}W_N(A_m)
	\right).
	\]
	
	\begin{theorem}\label{thm:rfif}
		Suppose that $W_i(x,y)=\bigl(L_i(x),F_i(x,y)\bigr),\,i=1,2,\ldots,N$,
		are contractive mappings satisfying $	L_i(x_{j(i)})=x_{i-1},
		L_i(x_{l(i)})=x_i,$
		and $	F_i(x_{j(i)},y_{j(i)})=y_{i-1},
		F_i(x_{l(i)},y_{l(i)})=y_i.$
		Then there exists a unique $A=(A_1,A_2,\ldots,A_N)\in \mathcal{H}(I\times \mathbb{R})^N$
		satisfying
		\[
		A_i
		=
		\bigcup_{m\in\mathcal{I}(i)}W_i(A_m),
		\qquad
		i=1,2,\ldots,N.
		\]
		Moreover, $G(f)
		=
		\bigcup_{i=1}^{N}A_i$
		is the graph of a continuous function $f:I\longrightarrow\mathbb{R}$
		satisfying $	f(x_i)=y_i,
		i=0,1,\ldots,N.$
	\end{theorem}
	
	\begin{definition}
		The continuous function $f$ obtained in
		Theorem~\ref{thm:rfif} is called the RFIF
		associated with the interpolation data $\Delta$ and the connection
		matrix $C$.
	\end{definition}
	
	\subsection{RFIF as a Fixed Point of the Read--Bajraktarevi\'c Operator}
	\label{subsec:rfif-rb}
	Denote by $C(I)$ the space of all real-valued continuous functions on
	$I$, endowed with the supremum norm $\|g\|_{\infty}
	=
	\max_{x\in I}|g(x)|.$
	Consider the closed subspace
	\[
	\mathcal{C}_{\Delta}(I)
	=
	\left\{
	g\in C(I):
	g(x_i)=y_i,\;
	i=0,1,\ldots,N
	\right\}.
	\]
	Define the Read--Bajraktarevi\'c operator $T:\mathcal{C}_{\Delta}(I)
	\longrightarrow
	\mathcal{C}_{\Delta}(I)$
	by
	\[
	(Tg)(x)
	=
	F_i
	\left(
	L_i^{-1}(x),
	g\left(L_i^{-1}(x)\right)
	\right),
	\qquad
	x\in I_i,
	\quad
	i=1,2,\ldots,N.
	\]
	\begin{theorem}\label{thm:rfif-rb}
		The operator $T:\mathcal{C}_{\Delta}(I)
		\longrightarrow
		\mathcal{C}_{\Delta}(I)$
		defined above is a contraction with contractivity factor $	s=\max_{1\leq i\leq N}s_i<1.$
		Its unique fixed point is the recurrent fractal interpolation function
		$f$ corresponding to the interpolation data $\Delta$. Consequently, $f$
		satisfies the self-referential equation
		\[
		f(L_i(x))
		=
		F_i(x,f(x)),
		\qquad
		x\in J_i,
		\quad
		i=1,2,\ldots,N.
		\]
	\end{theorem}
	
	The most widely used recurrent fractal interpolation functions are obtained
	by choosing $L_i(x)=a_i x+b_i$ and $F_i(x,y)=\alpha_i y+q_i(x),$ where $|\alpha_i|<1$ and $	q_i:J_i\longrightarrow\mathbb{R}$ is continuous function satisfying $q_i(x_{j(i)})
	=
	y_{i-1}-\alpha_i y_{j(i)} 
	\text{ ~and~  }
	q_i(x_{l(i)})
	=
	y_i-\alpha_i y_{l(i)}.$ The subinterval end-point conditions on $L_i$ yields,
	\[
	a_i
	=
	\frac{x_i-x_{i-1}}
	{x_{l(i)}-x_{j(i)}},
	\qquad
	b_i
	=
	\frac{
		x_{l(i)}x_{i-1}-x_{j(i)}x_i
	}{
		x_{l(i)}-x_{j(i)}
	}.
	\]
	The parameter $\alpha_i$ is referred to as the vertical scaling factor of
	the map $W_i$, and the vector $\boldsymbol{\alpha}
	=
	(\alpha_1,\alpha_2,\ldots,\alpha_N)
	\in(-1,1)^N$
	is called the scale vector of the recurrent IFS. If $q_i:J_i\longrightarrow\mathbb{R}, i=1,2,\ldots,N,$
	are affine maps, then the corresponding RFIF is termed an affine RFIF.
	In this case, $	q_i(x)=q_{i0}x+q_{i1},$
	where the recurrent join-up conditions $
	q_i(x_{j(i)})
	=
	y_{i-1}-\alpha_i y_{j(i)}
		\text{ ~and~  }
	q_i(x_{l(i)})
	=
	y_i-\alpha_i y_{l(i)}$
	yields
	\[
	q_{i0}
	=
	\frac{y_i-y_{i-1}}
	{x_{l(i)}-x_{j(i)}}
	-
	\alpha_i
	\frac{
		y_{l(i)}-y_{j(i)}
	}{
		x_{l(i)}-x_{j(i)}
	}
		\text{ ~and~  }
	q_{i1}
	=
	\frac{
		x_{l(i)}y_{i-1}-x_{j(i)}y_i
	}{
		x_{l(i)}-x_{j(i)}
	}
	-
	\alpha_i
	\frac{
		x_{l(i)}y_{j(i)}-x_{j(i)}y_{l(i)}
	}{
		x_{l(i)}-x_{j(i)}
	}.
	\]
	More generally, a useful choice of $q_i$ is $	q_i(x) = h\circ L_i(x)-\alpha_i b_i(x),  x\in J_i,$ where $h$ is "height function", which is a continuous interpolant to the prescribed data and
	$b_i\in C(J_i)$ is a base function satisfying $	b_i(x_{j(i)})=y_{j(i)},~
	b_i(x_{l(i)})=y_{l(i)}.$
    
	\section{Recurrent Fractal Functions of More General Nature}\label{sec3}
    
	In this section, we extend the construction of bounded fractal functions
	to the recurrent setting. By making suitable modifications to the
	construction of recurrent fractal interpolation functions discussed in
	the previous section, we remove both the continuity and interpolatory
	conditions inherent in an RFIF. Let $B(I)$ denote the linear space of all real-valued bounded functions defined on the compact interval $I$
	endowed with the norm $	\|g\|_{\infty}=	\sup\{|g(x)|:x\in I\}.$ Then $(B(I),\|\cdot\|_{\infty})$ is a Banach space.
	
    Let $x_0<x_1<\cdots<x_N$ be a partition of $I$. For each $i=1,2,\ldots,N$, let $	L_i:J_i\longrightarrow I_i$ be the contractive affine map $	L_i(x)=a_i x+b_i$ satisfying $L_i(x_{j(i)})=x_{i-1},~L_i(x_{l(i)})=x_i.$ For each $i=1,2,\ldots,N$, let $q_i:J_i\longrightarrow\mathbb{R}$ be a bounded function and define $	F_i:J_i\times\mathbb{R}\longrightarrow\mathbb{R}$ by 
    \begin{equation}\label{Fi-map}
        F_i(x,y)	=	\alpha_i y+q_i(x),
    \end{equation} where $	|\alpha_i|<1.$	In contrast to the construction of continuous RFIFs, we do not impose
	any join-up or interpolation conditions on the maps $F_i$. Now define $	W_i:J_i\times\mathbb{R}
	\longrightarrow
	I_i\times\mathbb{R}$ by
    \begin{equation}\label{Wi-map}
	    W_i(x,y)
	=
	\bigl(L_i(x),F_i(x,y)\bigr),
	\end{equation}
    For the definition of the Read--Bajraktarevi\'c operator, we regard the
	partition intervals as $I_i^-=[x_{i-1},x_i),	i=1,2,\ldots,N-1,$
	and $I_N=[x_{N-1},x_N],$ such that $	I=\bigcup_{i=1}^{N}I_i^-$ and the intervals are pairwise disjoint. Define the Read--Bajraktarevi\'c operator $T:B(I)\longrightarrow B(I)$
	by
	\[
	(Tg)(x)
	=
	F_i
	\left(
	L_i^{-1}(x),
	g\left(L_i^{-1}(x)\right)
	\right),
	\qquad
	x\in I_i^-.
	\]
	Equivalently,
\begin{equation}\label{Tg}
		(Tg)(x)
	=
	\alpha_i
	g\left(L_i^{-1}(x)\right)
	+
	q_i\left(L_i^{-1}(x)\right),
	\qquad
	x\in I_i^-.
\end{equation}
	\begin{theorem}\label{bounded-rff}
		The mapping $T:B(I)\longrightarrow B(I)$ defined by \eqref{Tg} is a contraction. Consequently, it possesses a
		unique fixed point $	f\in B(I),$ which satisfies the recurrent self-referential equation
	 \begin{equation}\label{eq3.2}
	 		f(L_i(x))
	 	=
	 	\alpha_i f(x)+q_i(x),
	 	\qquad
	 	x\in J_i,
	 	\quad
	 	i=1,2,\ldots,N.
	 \end{equation}
	\end{theorem}
	
	\begin{remark}
		Instead of constant vertical scaling factors, one may consider bounded
		scaling functions $\alpha_i:J_i\longrightarrow\mathbb{R}, i=1,2,\ldots,N.$ We may assume that each $\alpha_i$ is bounded and $\|\boldsymbol{\alpha}\|_{\infty}=
		\max_{1\leq i\leq N}
		\|\alpha_i\|_{\infty,J_i}
		<1,$
	
	\end{remark}
	\begin{remark}
		Similar to the construction of $\alpha$-fractal functions, let
		$g,b\in B(I)$ and define $q_i(x)
		=
		g(L_i(x))-\alpha_i b(x),x\in J_i,i=1,2,\ldots,N.$ Then the corresponding bounded recurrent fractal function satisfies
        \begin{equation*}
            f(L_i(x))
		=
		g(L_i(x))
		+
		\alpha_i(f-b)(x),
		\qquad x\in J_i.
        \end{equation*}
		The resulting function may be regarded as a recurrent
		$\alpha$-fractal analogue of $g$.
	\end{remark}

	\begin{theorem}\label{thm:rff-attractor}
		Let $f\in B(I)$ be the bounded recurrent fractal function obtained as the
		fixed point of the RB-operator in Theorem~\ref{bounded-rff}. For each
		$i=1,2,\ldots,N$, let $W_i$ defined in \eqref{Wi-map} and Suppose that $q_i:J_i\to\mathbb{R}$ is Lipschitz continuous with
		Lipschitz constant $Q_i$, and put $Q	=	\max_{1\leq i\leq N}Q_i,
		a=		\max_{1\leq i\leq N}a_i$ with $Q>0$, choose $\theta>0$ such that $0<\theta<\frac{1-a}{Q}.$ Then, the IFS $\{ J_i \times \mathbb{R}; W_i, i=1,\cdots, N\}$ is contractive with respect to a metric $d_\theta$ defined by .
		\[
		d_\theta
		\bigl((x,y),(x',y')\bigr)
		=
		|x-x'|+\theta|y-y'|.
		\]
		Moreover, the unique attractor of this IFS is the closure of the graph of $f$ .
	\end{theorem}

	\begin{remark}\label{rem:continuous-rff-graph}
		If the bounded recurrent fractal function $f$ in Theorem~\ref{thm:rff-attractor} is continuous on $I$, then its graph $G(f) = \{(x,f(x)):x\in I\}$ is closed. Hence, by Theorem~\ref{thm:rff-attractor}, for $i=1,2,\ldots,N; A_i = G_i(f). $ Consequently, \[ \bigcup_{i=1}^{N}A_i = \bigcup_{i=1}^{N}G_i(f) = G(f). \] Thus, in the continuous case, the attractor of the recurrent system coincides with the graph of the recurrent fractal function.
	\end{remark}
	The following result provides an error estimate for approximating a
	prescribed bounded function by a bounded recurrent fractal function.
	It may be regarded as a recurrent analogue of the standard
	collage-type estimate.
	\begin{theorem}\label{thm:recurrent-collage}
		Let $	\Phi\in B(I)$ be a bounded function and let $	T:B(I)\longrightarrow B(I)$ be the recurrent Read--Bajraktarevi\'c operator whose unique fixed point
		is the bounded recurrent fractal function $f$. Suppose that $\|\Phi-T\Phi\|_\infty<\varepsilon.$
		Then
		\[
		\|\Phi-f\|_\infty
		<
		\frac{\varepsilon}
		{1-|\alpha|_\infty},
		\]
		where $|\alpha|_\infty =	\max_{1\leq i\leq N}|\alpha_i|.$
	\end{theorem}
	As a prelude to the study of the discontinuities of bounded recurrent
	fractal functions, we recall the notion of oscillation.
	
	\begin{definition}\label{def:oscillation}
		Let $f:X\to\mathbb{R}$ be a function defined on a metric space
		$(X,d)$. For a nonempty set $U\subseteq X$, the oscillation of $f$
		on $U$ is defined by
		\[ \omega_f(U) = \sup_{x\in U}f(x) - \inf_{x\in U}f(x) = \sup_{x,y\in U}|f(x)-f(y)|.	\]
		The oscillation of $f$ at a point $x^*\in X$ is defined by
		\[
		\omega_f(x^*)
		=
		\lim_{\varepsilon\to0^+}
		\omega_f\bigl(B_\varepsilon(x^*)\bigr),
		\]
		where $B_\varepsilon(x^*) =\{x\in X:d(x,x^*)<\varepsilon\}$ denotes the open ball centered at $x^*$ with radius $\varepsilon$.
	\end{definition}
	
	\begin{remark}\label{rem:oscillation-continuity}
		A function $f:X\to\mathbb{R}$ is continuous at $x^*\in X$ if and only if $	\omega_f(x^*)=0.$
	\end{remark}
		\begin{theorem}\label{thm:rff-riemann}
		Let $f\in B(I)$ be the bounded recurrent fractal function obtained as
		the unique fixed point of the RB-operator $T$ in
		Theorem~\ref{bounded-rff}. Assume that, for each
		$i=1,2,\ldots,N$, the function $	q_i:J_i\longrightarrow\mathbb{R}$ is continuous. Then the set of points of discontinuity of
		$f$ is at most countable and, consequently, is a Lebesgue null set.
		In particular, $f$ is Riemann integrable on $I$.
	\end{theorem}

	\begin{lemma}\cite[Theorem~2(C. Neumann)]{yosida2012functional}\label{lemma}
	Let $A$ be a bounded linear operator on a Banach space. If $\|A\|<1,$ then $I-A$ is invertible. Moreover,
	\[(I-A)^{-1}=\sum_{n=0}^{\infty}A^n.	\]
\end{lemma}

	Now, we assume that $q_i$ is Lipschitz continuous with Lipschitz constant $Q_i$ and $Q=\max_{1\leq i \leq N} Q_i$
	\begin{theorem}\label{thm:recurrent-moments}
		Let $f$ be the bounded recurrent fractal function given in
		Theorem~\ref{bounded-rff}. For $	m\in\mathbb{N}\cup\{0\}$ and $i=1,2,\ldots,N$, define
		\[
		f_{i,m}
		=
		\int_{I_i}x^m f(x)\,dx,
		\]
		and set
		\[\mathbf f_m		=	
		\begin{pmatrix}	f_{1,m}\\
			f_{2,m}\\
			\vdots\\
			f_{N,m}
		\end{pmatrix}.\]
		For each $i$, define $Q_{i,m}
		=
		\int_{I_i}
		x^m q_i\left(L_i^{-1}(x)\right)\,dx,$ and let
		\[
		\mathbf Q_m
		=
		\begin{pmatrix}
			Q_{1,m}\\
			Q_{2,m}\\
			\vdots\\
			Q_{N,m}
		\end{pmatrix}.
		\]
		Define $D_m=\operatorname{diag}\left(\alpha_1a_1^{m+1},\ldots,\alpha_Na_N^{m+1}\right),$ and, for $0\leq k<m$, $D_{m,k}
		=
		\operatorname{diag}
		\left(
		\alpha_1a_1^{k+1}b_1^{m-k},
		\ldots,
		\alpha_Na_N^{k+1}b_N^{m-k}
		\right).$ Then $I_N-D_mC$ is non singular and
		\[
			\mathbf f_m
			=
			\left(I_N-D_mC\right)^{-1}
			\left[
			\sum_{k=0}^{m-1}
			\binom{m}{k}
			D_{m,k}C\mathbf f_k
			+
			\mathbf Q_m
			\right].
		\]
	   In particular, $	\mathbf f_0
	   =
	   \left(I_N-D_0C\right)^{-1}
	   \mathbf Q_0.$ Moreover, the $m$th moment of $f$ on $I$ is $f_m
		=
		\int_Ix^mf(x)\,dx
		=
		\sum_{i=1}^{N}f_{i,m}.$
	\end{theorem}

    We next consider another consequence of the recurrent functional equation, namely the integral transforms of bounded recurrent fractal functions. Let $f:I\to\mathbb{R}$ be a bounded recurrent fractal function. Extend $f$ to $\mathbb{R}$ by setting $f(x)=0,~x\in\mathbb{R}\setminus I,$ and, for simplicity, denote the extended function again by $f$. Then, for a suitable kernel $K(x,s)$,
    \[\widehat{f}(s)=\int_{\mathbb{R}} K(x,s)f(x)\,dx=\int_I K(x,s)f(x)\,dx,\]
    since $f$ vanishes outside $I$.
	\section{Recurrent Fractal Histopolation}\label{sec4}

    	Let $x_0<x_1<\cdots<x_N $ be a sequence of strictly increasing knots and let $F=\{f_1,f_2,\ldots,f_N\}$ be a prescribed histogram, where $f_i\in\mathbb{R}$ represents the frequency corresponding to the class interval $I_i=[x_{i-1},x_i),~ i=1,2,\ldots,N-1, $ and $I_N=[x_{N-1},x_N].$  For each $i=1,2,\ldots,N$, let $h_i=x_i-x_{i-1}$ denote the corresponding step size.

	In recurrent fractal histopolation, instead of prescribing the values
	of a fractal function at the knots, we prescribe its average values
	over the basic intervals. More precisely, we seek a bounded integrable
	recurrent fractal function $f\in B(I)$ that satisfies the area-matching
	conditions
	\begin{equation}\label{areacond}
		\int_{x_{i-1}}^{x_i}f(x)\,dx = h_i f_i, \qquad i=1,2,\ldots,N.
	\end{equation}
	Such a function will be called a recurrent fractal histopolant
	corresponding to the histogram $F$. Consider $L_i$ and $F_i$
	\begin{equation}\label{maps}
		L_i(x)=a_i x+b_i, \quad F_i(x,y)=\alpha_i y+q_i(x), \quad x\in J_i, \quad i=1,2,\ldots,N.
	\end{equation}
    where $	0<a_i<1, ~|\alpha_i|<1,$ and $q_i:J_i\longrightarrow\mathbb{R}$ is Lipschitz continuous. From Theorem \ref{bounded-rff} and Theorem \ref{thm:rff-riemann} it follows that corresponding fractal function is Riemann integrable and satisfies
	\[f(x) = \alpha_i f\left(L_i^{-1}(x)\right) + q_i\left(L_i^{-1}(x)\right), \qquad x\in I_i,\]

    	\begin{proposition}\label{prop:recurrent-histopolation}
		Let $x_0<x_1<\cdots<x_N $ be a sequence of strictly increasing knots and let $F=\{f_1,f_2,\ldots,f_N\}$ be a prescribed histogram,  Consider the maps given in \eqref{maps}. Then the corresponding bounded recurrent fractal function $f$ solves
		the histopolation problem give in \eqref{areacond} if and only if function satisfies
	\begin{equation}\label{condq_i}
      	\int_{J_i}q_i(t)\,dt = \frac{ h_i f_i - a_i\alpha_i \displaystyle\sum_{r=1}^{N} C_{ir}h_r f_r }{a_i}, \qquad i=1,2,\ldots,N.
      \end{equation}
	\end{proposition}
	 Next suppose that we are interested in constructing a continuous recurrent fractal histopolant corresponding to the strictly increasing knots $x_0<x_1<\cdots<x_N $ and the histogram $F=\{f_1,f_2,\ldots,f_N\}.$ For $r=0,1,\ldots,N,$ let $y_r=f(x_r)$ denote the values of the desired continuous histopolant at the knots. Recall from Section \ref{sec2} that the recurrent fractal function corresponding to the map \eqref{maps} is continuous if 
	  \begin{equation}\label{qi}
	  	\begin{cases} q_i(x_{j(i)}) = y_{i-1}-\alpha_i y_{j(i)},\\[2mm] q_i(x_{l(i)}) = y_i-\alpha_i y_{l(i)}, 
	  	\end{cases} \qquad i=1,2,\ldots,N. 
	  \end{equation}
      Now we are ready to provide the following proposition.
      
	 \begin{proposition}\label{prop:continuous-recurrent-histopolation} Let $ x_0<x_1<\cdots<x_N $ be a sequence of strictly increasing knots and let $ F=\{f_1,f_2,\ldots,f_N\} $ be a prescribed histogram. Consider the maps given in \eqref{maps}. Assume that scaling factor are fixed i.e., $| \alpha_i| <1$. The corresponding fractal function is continuous and solves the histopolation problem in \eqref{areacond}, if the function $q_i$ satisfies the system of equation governed by \eqref{condq_i}  and \eqref{qi} 
	 \end{proposition}

\section*{Declaration}\noindent
	\textbf{Conflicts of interest.} We do not have any conflict of interest.\\
	\\
	\noindent
	\textbf{Data availability:} No data were used to support this study.\\
	\\
	\noindent
	\textbf{Code availability:} Not applicable\\
	\\
	\noindent
	\textbf{Authors' contributions:} Each author contributed equally to this manuscript.

	\bibliographystyle{elsarticle-num}
	\bibliography{Reference}
    \end{document}